\documentclass[11pt]{article}
\usepackage{graphicx,mathrsfs,amssymb}
\usepackage{amsmath,latexsym, color, amsbsy,amssymb}
\usepackage{amsthm}
\theoremstyle{plain}
\usepackage{graphicx}

\makeatletter
\newcommand{\superimpose}[2]{%
	{\ooalign{$#1\@firstoftwo#2$\cr\hfil$#1\@secondoftwo#2$\hfil\cr}}}
\makeatother

\def\ds{\displaystyle}
\def\forall{\hbox{for all}~}

\def\n{\noindent}

\def\R{\mathbb{R}}

\def\vs{\vskip 2em}

\def\v{\vskip 1em}

\def\bega{\begin{array}}
	\def\enda{\end{array}}
\def\begi{\begin{itemize}}
	\def\endi{\end{itemize}}

\def\bel{\begin{equation}\label}
	\def\eeq{\end{equation}}
\def\sqr#1#2{\vbox{\hrule height .#2pt
		\hbox{\vrule width .#2pt height #1pt \kern #1pt
			\vrule width .#2pt}\hrule height .#2pt }}

\def\square{\sqr74}
\def\endproof{\hphantom{MM}\hfill\llap{$\square$}\goodbreak}

\newtheorem*{theorem*}{Theorem}

\begin{document}
	\title{\bf  Hausdorff measure of zeros of polynomials}\vs
	\author{\it Andrew Murdza, Khai T. Nguyen, and  Etienne Phillips\\
\\
		{\small Department of Mathematics, North Carolina State University}\\
		{\small e-mails: ~apmurdza@ncsu.edu, ~khai@math.ncsu.edu,  ~ecphill6@ncsu.edu}
}
\maketitle
\begin{abstract} 
The paper  provides an elementary proof  establishing a sharp universal bound on the $(d-1)$-Hausdorff measure of the zeros of any nontrivial multivariable polynomial $p:\mathbb{R}^d\to\mathbb{R}$ within a $d$-dimensional cube of size $r$. This bound depends solely on the parameter $r$, the dimension $d$, and the degrees of $p$.
\quad\\
		\quad\\
		{\footnotesize
		{\bf Keywords.}  Hausdorff measure,  zeros of polynomials
		\medskip
		
		\n {\bf AMS Mathematics Subject Classification.} 12D10, 28A78
		}
	\end{abstract}
%
	
The fundamental theorem of algebra, also known as d'Alembert's theorem or the d'Alembert-Gauss theorem, asserts that the number of roots (counting multiplicities) for a nontrivial polynomial of a single variable is equal to its degree. Various proofs have been developed in attempts to comprehend the theorem from different mathematical perspectives, including Topology (refer to \cite{A} or \cite{ST}), Real and Complex Analysis (refer to \cite{TK}), and Algebra (refer to \cite{H} or \cite{Z}). These diverse approaches create non-trivial intersections with multiple branches of mathematics. A compilation of proofs can be found in \cite{FR}. The  focus of this paper is to address the following natural question:
\medskip

\qquad {\bf Question.} {\it Can an analogous result to the fundamental theorem of algebra be provided for multivariable polynomials?}
\medskip

For every nontrivial polynomial $p:\R^d\to\R$, the zero set of  $p$, denoted by,
\[
\mathcal{Z}_{p}~\doteq~\left\{x=(x_1,\cdots,x_d)\in\R^d: p(x)~=~0\right\},
\]
has Hausdorff dimension less than or equal $d-1$.  This claim remains valid for every non-identically zero analytic function and certain quasi-analytic functions as outlined in \cite{GS}. Multiple elementary proofs, detailed in \cite[Remark 5.23]{PK}, include applications of Fubini's theorem as demonstrated in \cite[Corollary 5.2]{BS}, induction with respect to dimension in \cite[Section 4.1]{KP}, and coverage of the set with a countable set of smooth submanifolds using the implicit function theorem \cite{M}. Moreover, it can be deduced from the Weierstrass Preparatory Theorem \cite{GR}, Lojasiewicz's result on the existence of Whitney stratification of real analytic sets \cite{L}, or the desingularization theorems by Bierstone and Milman \cite{BM}. Despite these approaches, the size of the zero set $\mathcal{Z}_p$ in the analogy of the fundamental theorem of algebra remains unaddressed. Here, we establish a sharp universal bound on the $(d-1)$-Hausdorff measure of $\mathcal{Z}_p$ within $d$-dimensional cubes.\medskip 

\begin{theorem*}
For every nontrivial polynomial $p:\R^d\to\R$, it holds 
\bel{H}
\mathcal{H}^{d-1}\big(\mathcal{Z}_p\cap [a,b]^d\big)~\leq~\left(\sum_{k=1}^{d}\deg_{x_k}p\right)\cdot (b-a)^{d-1},\qquad a<b,
\eeq
where $\deg_{x_k}p$ represents the degree of $p$ in the variable $x_k$.
\end{theorem*}
{\bf Proof.}
{\bf 1.} By the fundamental theorem of algebra, (\ref{H}) holds for $d=1$. Thus, we only need to prove (\ref{H}) for $d\geq 2$. Given a nontrivial polynomial $p:\R^{d}\to\R$, we set $\kappa_d\doteq \deg_{x_d}p$ and write 
\[
p(x)~=~q_0(x_1,\cdots,x_{d-1})+q_1(x_1,\cdots,x_{d-1})x_{d}+\cdots+ q_{\kappa_{d}}(x_1,\cdots,x_{d-1})x^{\kappa_{d}}_{d}.
\]
Consider the $\kappa_d$ disjoint subsets $E_{\ell}$ of $\mathcal{Z}_p$ for $1\leq \ell\leq \kappa_d$, defined as follows: 
\[
E_{\ell}~=~\left\{x\in [a,b]^d:{\partial^{\ell}\over \partial x^{\ell}_d}p(x)\neq 0, {\partial^{j}\over \partial x^j_d}p(x)=0~~\forall 0\leq j\leq \ell-1\right\}.
\]
For every $k\in\{1,\cdots, d\}$, let $\pi_k:\R^{d}\to\R^{d-1}$ be the coordinate projection such that 
\bel{pi-k}
\pi_k(x)~=~(x_1,\cdots x_{k-1}, x_{k+1},\cdots, x_{d}),\qquad\forall x\in \R^{d},
\eeq
and 
\bel{T-k}
\mathcal{T}_k~=~\left\{y\in \R^{d-1}:p(x)=0~~\forall x\in\pi_k^{-1}(y)\right\}.
\eeq
Recalling that $\mathcal{Z}_{q_{\kappa_{d}}}$ is the zero set of the polynomial $q_{\kappa_d}$, we decompose $\mathcal{Z}_{p}\cap [a,b]^d$ into two subsets by
\bel{P-Z}
\mathcal{Z}_{p}\cap [a,b]^d~=~ \mathcal{E}_d\bigcup \left\{x\in \mathcal{Z}_{p}\cap [a,b]^d:\pi_d(x)\in \mathcal{Z}_{q_{\kappa_{d}}} \right\},\qquad \mathcal{E}_d~\doteq~\bigcup_{\ell=1}^{\kappa_d}E_{\ell}.
\eeq
For every $x\in E_{\ell}$,  by the implicit function theorem, there exist an open neighborhood $U_x$ of $\pi_d(x)$ in $\R^{d-1}$  and a smooth map $\varphi_x:U_{x}\to \R$ satisfying $\ds{\partial^{\ell-1}\over\partial x_d^{\ell-1} }p(y,\varphi_x(y))=0$ for all $y\in U_x$. For every small neighborhood $V_x$ of $x$ in $\R^d$ such that $V_x\subset U_x\times\varphi_x(U_x)$,  using  the area formula (see e.g. in \cite[Theorem 2.71]{AM}), we estimate   
\bel{edq}
\begin{split}
\mathcal{H}^{d-1}\big(E_{\ell}\cap V_x\big)&~=~\int_{\pi_{d}(E_{\ell}\cap V_x)}\sqrt{1+|\nabla \varphi_x(y)|^2}~dy~\leq~\int_{\pi_{d}(E_{\ell}\cap V_x)}1+\sum_{k=1}^{d-1}\left|{\partial\over \partial y_k}\varphi_x(y)\right|~dy\\
&~=~\int_{ [a,b]^{d-1}}\chi_{\pi_{d}(E_{\ell}\cap V_x)}(y)dy+\sum_{k=1}^{d-1}\int_{\pi_{d}(E_{\ell}\cap V_x)}\left|{\partial\over \partial y_k}\varphi_x(y)\right|~dy\\
&~=~\int_{ [a,b]^{d-1}}\mathcal{H}^0\left(\pi_d^{-1}(y)\cap E_{\ell}\cap V_x\right)dy+\sum_{k=1}^{d-1}\int_{\pi_{d}(E_{\ell}\cap V_x)}\left|{\partial\over \partial y_k}\varphi_x(y)\right|~dy.
\end{split}
\eeq
For every $k\in \{1,\cdots, d-1\}$ and $y\in \R^{d-1}$, writing   $y^{-}_{k}=(y_1,\cdots y_{k-1}, y_{k+1},\cdots, y_{d-1})$ and $I_{V_x}(y^-_k)=\left\{t\in [a,b]:(y_1,\cdots y_{k-1}, t,y_{k+1},\cdots, y_{d-1})\in \pi_{d}(E_{\ell}\cap V_x)\right\}$, by Fubini's theorem and the area formula, we obtain 
\bel{et1}
\begin{split}
\int_{\pi_{d}(E_{\ell}\cap V_x)}\left|{\partial\over \partial y_k}\varphi_x(y)\right|~dy&~=~\int_{[a,b]^{d-2}}\int_{I_{V_x}(y^-_k)}\left|{\partial\over \partial y_k}\varphi_x(y)\right|dy_{k}dy^{-}_{k}\\
&~=~\int_{[a,b]^{d-2}}\int_{a}^{b}\mathcal{H}^{0}\left(\left\{y_k\in I_{V_x}(y^-_k):\varphi_x(y)=\tau\right\}\right)d\tau dy^{-}_{k}\\
&~=~\int_{[a,b]^{d-2}}\int_{a}^{b}\mathcal{H}^{0}\left(\pi_k^{-1}(y_k^-,\tau)\cap \big[E_{\ell}\cap V_x\big]\right)d\tau dy^{-}_{k}.
\end{split}
\eeq
Notice that if $(y_k^-,\tau_0)\in\mathcal{T}_{k}\cap \pi_k\left(E_{\ell}\cap V_x\right) $ then 
\[
\begin{cases}
\varphi_x(y_1,\cdots y_{k-1}, t,y_{k+1},\cdots, y_{d-1})~=~\tau_0,\qquad\forall t\in I_{V_x}(y^-_k),\\[3mm]
p(y_1,\cdots y_{k-1}, t,y_{k+1},\cdots, y_{d-1},\tau_0)~=~0,\qquad\forall t\in \R.
\end{cases}
\]
In this case, we have 
\[
\mathcal{H}^{0}\left(\pi_k^{-1}(y_k^-,\tau)\cap \big[E_{\ell}\cap V_x\big]\right)~=~0,\qquad\forall \tau\neq \tau_0,
\]
and this implies  
\[
\int_{a}^{b}\mathcal{H}^{0}\left(\pi_k^{-1}(y_k^-,\tau)\cap \big[E_{\ell}\cap V_x\big]\right)d\tau~=~0.
\]
In particular, from (\ref{et1}), we derive
\[
\begin{split}
\int_{\pi_{d}(E_{\ell}\cap V_x)}\left|{\partial\over \partial y_k}\varphi_x(y)\right|~dy&~\leq~\int_{[a,b]^{d-1}\backslash[\mathcal{T}_k\cap \pi_k\left(E_{\ell}\cap V_x\right) ]}\mathcal{H}^0\left(\pi_k^{-1}(z)\cap (E_{\ell}\cap V_x)\right)dz\\
&~=~\int_{[a,b]^{d-1}\backslash\mathcal{T}_k}\mathcal{H}^0\left(\pi_k^{-1}(z)\cap (E_{\ell}\cap V_x)\right)dz,
\end{split}
\]
%
and (\ref{edq}) yields 
\[
\mathcal{H}^{d-1}\big(E_{\ell}\big)~\leq~\int_{ [a,b]^{d-1}}\mathcal{H}^0\left(\pi_d^{-1}(y)\cap E_{\ell}\right)dy+\sum_{k=1}^{d-1}\int_{[a,b]^{d-1}\backslash\mathcal{T}_k}\mathcal{H}^0\left(\pi_k^{-1}(z)\cap E_{\ell}\right)dz.
\]
Thus, recalling (\ref{P-Z}), we get 
\bel{H-E-l}
\begin{split}
\mathcal{H}^{d-1}\big(\mathcal{E}_d \big)&~=~\sum_{\ell=1}^{\kappa_d}\mathcal{H}^{d-1}\left(E_{\ell}\right)~\leq~\sum_{\ell=1}^{\kappa_d}\sum_{k=1}^{d}\int_{[a,b]^{d-1}\backslash\mathcal{T}_k}\mathcal{H}^0\left(\pi_k^{-1}(y)\cap E_{\ell}\right)dy\\
&~=~\sum_{k=1}^{d}\int_{[a,b]^{d-1}\backslash\mathcal{T}_k}\mathcal{H}^0\left(\pi_k^{-1}(y)\cap \mathcal{E}_{d}\right)dy.
\end{split}
\eeq

{\bf  2.} Let's now estimate $\mathcal{H}^{d-1}\left(\left\{x\in \mathcal{Z}_{p}\cap [a,b]^d:\pi_d(x)\in \mathcal{Z}_{q_{\kappa_{d}}}\right\}\right)$. For every $\ell\in \{0,\cdots,\kappa_d-1\}$, we define 
\[
\Omega^{\ell}_{q_{\kappa_{d}}}~\doteq~\left\{y\in \R^{d-1}:q_{\ell}(y)\neq 0~~\mathrm{and}~~q_{j}(y)=0~~\forall \ell+1\leq j\leq \kappa_d\right\},
\]
and 
\[
\Omega^{c}_{q_{\kappa_d}}~\doteq~\left(\mathcal{Z}_{q_{\kappa_d}}\cap [a,b]^{d-1}\right)\backslash\bigcup_{\ell=1}^{\kappa_d-1}\Omega^{\ell}_{q_{\kappa_{d}}}~=~\left\{y\in [a,b]^{d-1}:q_{j}(y)=0~~\forall 0\leq j\leq \kappa_d\right\}.
\]
By  the continuous dependence of roots on the coefficients (see e.g.~in \cite{NP}), there exist $\ell$ continuous functions $\gamma_i:  \Omega^{\ell}_{q_{\kappa_{d}}}\to \mathbb{C}$ such that  
\[
\left\{x\in \mathcal{Z}_{p}\cap [a,b]^d:\pi_d(x)\in \Omega^{\ell}_{q_{\kappa_{d}}}\right\}~\subseteq~\bigcup_{i=1}^{\ell}\left\{(y,\gamma_{i}(y)):y\in \Omega^{\ell}_{q_{\kappa_{d}}}\right\}~\subseteq~\mathcal{Z}_{p}.
\]
Since $q_{\kappa_d}$ is a nontrivial polynomial, one has  
\[
\mathcal{L}^{d-1}\big(\Omega^{\ell}_{q_{\kappa_{d}}}\big)~\leq~ \mathcal{L}^{d-1}\big( {\mathcal{Z}}_{q_{\kappa_{d}}}\big)~=~ 0,
\]
and  the continuity of $\gamma_i$ implies that 
\[
\mathcal{H}^{d-1}\left(\left\{x\in \mathcal{Z}_{p}\cap [a,b]^d:\pi_d(x)\in \Omega^{\ell}_{q_{\kappa_{d}}}\right\}\right)~\leq~\sum_{i=1}^{\ell}\mathcal{H}^{d-1}\left(\left\{(y,\gamma_{i}(y)):y\in \Omega^{\ell}_{q_{\kappa_{d}}}\right\}\right)~=~0.
\]
In particular, we obtain
\[
\begin{split}
\mathcal{H}^{d-1}\left(\left\{x\in \mathcal{Z}_{p}\cap [a,b]^d:\pi_d(x)\in \mathcal{Z}_{q_{\kappa_{d}}}\right\}\right)&=~\mathcal{H}^{d-1}\left(\left\{x\in \mathcal{Z}_{p}\cap [a,b]^d:\pi_d(x)\in \Omega^{c}_{q_{\kappa_d}}\right\}\right)\\
&~=~\mathcal{H}^{d-1}\left(\Omega^{c}_{q_{\kappa_d}}\times [a,b]\right),
\end{split}
\]
and (\ref{P-Z}), (\ref{H-E-l}) yield
\bel{ett2}
\mathcal{H}^{d-1}\left(\mathcal{Z}_{p}\cap [a,b]^d\right)~\leq~\mathcal{H}^{d-1}\left(\Omega^{c}_{q_{\kappa_d}}\times [a,b]\right)+\sum_{k=1}^{d}\int_{[a,b]^{d-1}\backslash\mathcal{T}_k}\mathcal{H}^0\left(\pi_k^{-1}(y)\cap \mathcal{E}_{d}\right)dy.
\eeq
To complete this step, we  notice that $\Omega^{c}_{q_{\kappa_d}}\times [a,b]$  and $\mathcal{E}_d$ are subsets of $\mathcal{Z}_p$ and have an empty intersection. 
\medskip

{\bf  3.} By the same argument in the previous steps  for the nontrivial polynomial $p_{d-1}\doteq q_{\kappa_d}:\R^{d-1}\to\R$, there exist a nontrivial polynomial $p_{d-2}:\R^{d-2}\to\R$ and two subsets $\mathcal{E}_{d-1}\subseteq \Omega^{c}_{p-1}\times [a,b]$ and $\Omega^{c}_{p_{d-2}}\subseteq [a,b]^{d-2}\cap \mathcal{Z}_{p_{d-2}}$ such that $\Omega^{c}_{p_{d-2}}\times [a,b]^2$ and $\mathcal{E}_{d-1}$  are subsets of $\mathcal{Z}_p$ and have an empty intersection, and 
\[
\mathcal{H}^{d-1}\left(\Omega^{c}_{p_{d-1}}\times [a,b]\right)~\leq~\mathcal{H}^{d-1}\left(\Omega^{c}_{p_{d-2}}\times [a,b]^2\right)+\sum_{k=1}^{d}\int_{[a,b]^{d-1}\backslash\mathcal{T}_k}\mathcal{H}^0\left(\pi_k^{-1}(y)\cap \mathcal{E}_{d-1}\right)dy.
\]
Thus, (\ref{ett2}) implies
\[
\mathcal{H}^{d-1}\left(\mathcal{Z}_{p}\cap [a,b]^d\right)~\leq~\mathcal{H}^{d-1}\left(\Omega^{c}_{p_{d-2}}\times [a,b]^2\right)+\sum_{k=1}^{d}\sum_{i=0}^{1}\int_{[a,b]^{d-1}\backslash\mathcal{T}_k}\mathcal{H}^0\left(\pi_k^{-1}(y)\cap \mathcal{E}_{d-i}\right)dy.
\]
Continuing this procedure, we can construct $d$ disjoint subsets $\mathcal{E}_{i}$ of $\mathcal{Z}_p\cap [a,b]^{d}$ with $i\in {0,1,\cdots, d-1}$ such that
\bel{ettf}
\begin{split}
\mathcal{H}^{d-1}\left(\mathcal{Z}_{p}\cap [a,b]^d\right)&~\leq~\sum_{k=1}^{d}\sum_{i=0}^{d-1}\int_{[a,b]^{d-1}\backslash\mathcal{T}_k}\mathcal{H}^0\left(\pi_k^{-1}(y)\cap \mathcal{E}_{d-i}\right)dy\\
&~=~\sum_{k=1}^{d}\int_{[a,b]^{d-1}\backslash\mathcal{T}_k}\mathcal{H}^0\left(\pi_k^{-1}(y)\cap \left(\bigcup_{i=0}^{d-1}\mathcal{E}_{d-i}\right)\right)dy\\
&~\leq~\sum_{k=1}^{d}\int_{[a,b]^{d-1}\backslash\mathcal{T}_k}\mathcal{H}^0\left(\pi_k^{-1}(y)\cap\mathcal{Z}_p\right)dy.
\end{split}
\eeq
Finally, by the fundamental theorem of algebra,  one has  
\[
\mathcal{H}^0\left(\pi_k^{-1}(y)\cap\mathcal{Z}_p\right)~\leq~\deg_{x_k}p,\qquad k\in \{1,\cdots, d\}, y\in [a,b]^{d-1}\backslash\mathcal{T}_k,
\]
and (\ref{ettf}) yields 
\[
\mathcal{H}^{d-1}\left(\mathcal{Z}_{p}\cap [a,b]^d\right)~\leq~\sum_{k=1}^{d}\int_{[a,b]^{d-1}\backslash\mathcal{T}_k}\deg_{x_k}p~dy~\leq~\left(\sum_{k=1}^{d}\deg_{x_k}p\right)\cdot (b-a)^{d-1}.
\]
The proof is complete.
\endproof
{\bf Remark.} {\it The constant $\ds\sum_{k=1}^{d}\deg_{x_k}p$ in (\ref{H}) is sharp. Indeed, consider the sequence of polynomials $p_n:\R^{d}\to \R$ such that $$p_n(x)~=~x_{1}x_2\cdots x_{d}-1/n,\qquad x\in \R^d,n\geq 1.
$$
We have that $\ds\sum_{k=1}^{d}\deg_{x_k}p_n=d$ for all $n\geq 1$ and $\ds\lim_{n\to+\infty}\mathcal{H}^{d-1}\left(\mathcal{Z}_{p_n}\bigcap [0,1]^d\right)=d$.
}
\v

{\small {\bf Acknowledgments.}  This research  was partially supported by National Science Foundation grant DMS-2154201. The authors would like to express their sincere gratitude to  Prof. J\'anos Koll\'ar for his insightful suggestions. }

%
%
%

\end{document}